\documentclass[12pt]{amsart}
\usepackage{latexsym,amssymb,amsthm,upref,mathrsfs}
\usepackage{amsmath}
\usepackage{cite}
\usepackage[top=0.8in, left=1in, right=1in, bottom=1in]{geometry}
\usepackage{listings}

\lstdefinelanguage{GAP}{%
  morekeywords={%
    Assert,Info,IsBound,QUIT,%
    TryNextMethod,Unbind,and,break,%
    continue,do,elif,%
    else,end,false,fi,for,%
    function,if,in,local,global,%
    mod,not,od,or,%
    quit,rec,repeat,return,%
    then,true,until,while,proc,options%
  },%
  sensitive,%
  morecomment=[l]\#,%
  morestring=[b]",%
  morestring=[b]',%
}[keywords,comments,strings]

\newcommand{\cl}{C \kern -0.1em \ell}

\makeatletter
\newcommand\GL{\mathop{\operator@font GL}\nolimits}
\makeatother

\makeatletter
\newcommand\Mat{\mathop{\operator@font Mat}\nolimits}
\makeatother

\makeatletter
\newcommand\sgn{\mathop{\operator@font sgn}\nolimits}
\makeatother

\makeatletter
\newcommand\spn{\mathop{\operator@font span}\nolimits}
\makeatother

\makeatletter
\newcommand\tr{\mathop{\operator@font tr}\nolimits}
\makeatother

\newcommand{\BR}{\mathbb{R}}
\newcommand{\BC}{\mathbb{C}}
\newcommand{\BH}{\mathbb{H}}

\newcommand{\BZ}{\mathbb{Z}}

\newcommand{\BK}{\mathbb{K}}

\newcommand{\be}{\mathbf{e}}

\newcommand{\ed}{\end{document}}

\newcommand{\cb}[1]{\mathcal{#1}}

\newcommand{\clpq}[2]{\cl_{#1,#2}}

\newcommand{\fpower}[1]{{}^2 \kern -0.2em #1} %left power

\newcommand{\iu}{{\underline{i}}}

\newcommand{\End}{\mathrm{End}}

\theoremstyle{plain}
\newtheorem{theorem}{Theorem}
\newtheorem{corollary}{Corollary}
\newtheorem{lemma}{Lemma}
\newtheorem{proposition}{Proposition}

\newtheorem*{StructureTheorem}{Structure Theorem}
\theoremstyle{definition}
\newtheorem{definition}{Definition}

\begin{document}
\vspace*{-4ex}
\title{On the Structure Theorem of Clifford Algebras}

\maketitle
\begin{center}
\author{{\bf Rafa\l \ Ab\l amowicz}\,$^a$}\\
\vspace{15pt}
\small
\thispagestyle{empty}
$^a$ Department of Mathematics, Tennessee Technological University \\
Cookeville, TN 38505, U.S.A. \\
\texttt{rablamowicz@tntech.edu}, \texttt{http://math.tntech.edu/rafal/}
\end{center}

\noindent
\begin{abstract}
In this paper, theory and construction of spinor representations of real Clifford algebras $\clpq{p}{q}$ in minimal left ideals are reviewed. Connection with a general theory of semisimple rings is shown. The actual computations can be found in, for example, \cite{spinorrepr}.
\end{abstract}

\noindent
{\small

}
\vspace*{-4ex}
\tableofcontents

\section{Introduction}
\label{sect:sect1}
Theory of spinor representations of real Clifford algebras $\clpq{p}{q}$ over a quadratic space $(V,Q)$ with a nondegenerate quadratic form $Q$ of signature $(p,q)$ is well known~\cite{chevalley,lam,lounesto,porteous}. The purpose of this paper is to review the structure theorem of these algebras in the context of a general theory of semisimple rings culminating with Wedderburn-Artin Theorem~\cite{rotman}. 

Section~2 is devoted to a short review of general background material on the theory of semisimple rings and modules as a generalization of the representation theory of group algebras of finite groups \cite{james, rotman}. While it is well-known that Clifford algebras $\clpq{p}{q}$ are associative finite-dimensional unital semisimple $\BR$-algebras, hence the representation theory of semisimple rings \cite[Chapter 7]{rotman} applies to them, it is also possible to view these algebras as twisted group algebras $\BR^t[(\BZ_2)^n]$ of a finite group $(\BZ_2)^n$ \cite{ablamowicz1,ablamowicz2,ablamowicz3,albuquerque,downs,majid}. While this last approach is not pursued here, for a connection between Clifford algebras $\clpq{p}{q}$ and finite groups, see 
\cite{ablamowicz2,ablamowicz3,brasov,brown,maduranga,maduranga2,walley} and references therein. 

In Section~3, we state the main Structure Theorem on Clifford algebras $\clpq{p}{q}$ and relate it to the general theory of semisimple rings, especially to the Wedderburn-Artin theorem. For details of computation of spinor representations, we refer to \cite{spinorrepr} where these computations were done in great detail by hand and by using \texttt{CLIFFORD}, a Maple package specifically designed for computing and storing spinor representations of Clifford algebras $\clpq{p}{q}$ for $n=p+q\leq 9$ \cite{clifford,clifford2}.

Our standard references on the theory of modules, semisimple rings and their representation is 
\cite{rotman}; for Clifford algebras we use~\cite{chevalley,lam,lounesto} and references therein; on representation theory of finite groups we refer to \cite{james,passman} and for the group theory we 
refer to~\cite{dornhoff, gorenstein, mckay, rotman}.

\section{Introduction to Semisimple Rings and Modules}
\label{sect:sect2}

This brief introduction to the theory of semisimple rings is based on \cite[Chapter 7]{rotman} and it is stated in the language of left $R$-modules. Here, $R$ denotes an associative ring with unity~$1$. We omit proofs as they can be found in Rotman~\cite{rotman}. 

\begin{definition}
Let $R$ be a ring. A \textbf{left $R$-module} is an additive abelian group $M$ equipped with \textbf{scalar multiplication} $R \times M \rightarrow M$, denoted $(r,m) \mapsto rm,$ such that the following axioms hold for all $m,m' \in M$ and all $r,r' \in R:$
\begin{itemize}
\item[(i)] $r(m+m')=rm+rm',$
\item[(ii)] $(r+r')m=rm+r'm,$
\item[(iii)] $(rr')m=r(r'm),$
\item[(iv)] $1m=m.$
\end{itemize}
Left $R$-modules are often denoted by $_RM.$
\end{definition}

In a similar manner one can define a \textbf{right $R$-module} with the action by the ring elements on~$M$ from the right. When $R$ and $S$ are rings and $M$ is an abelian group, then $M$ is a \textbf{$(R,S)$-bimodule}, denoted by $_RM_S$, if $M$ is a left $R$-module, a right $S$-module, and the two scalar multiplications are related by an associative law: $r(ms) = (rm)s$ for all $r \in R, m \in M,$ and $s \in S.$

We recall that a spinor left ideal $S$ in a simple Clifford algebra $\clpq{p}{q}$ by definition carries an irreducible and faithful representation of the algebra, and it is defined as $\clpq{p}{q}f$ where $f$ is a primitive idempotent in $\clpq{p}{q}$. Thus, as it is known from the Structure Theorem (see Section~\ref{sect:sect3}), that these ideals are $(R,S)$-bimodules where $R=\clpq{p}{q}$ and $S=f\clpq{p}{q}f$. Similarly, the right spinor modules $f\clpq{p}{q}$ are $(S,R)$-bimodules. Notice that the associative law mentioned above is automatically satisfied because $\clpq{p}{q}$ is associative.

We just recall that when $k$ is a field, every finite-dimensional $k$-algebra $A$ is both 
\textbf{left} and \textbf{right noetherian}, that is, any ascending chain of left and right ideals stops (the \textbf{ACC ascending chain condition}). This is important for Clifford algebras because, eventually, we will see that every Clifford algebra can be decomposed into a finite direct sum of left spinor $\clpq{p}{q}$-modules 
(ideals). For completeness we mention that every finite-dimensional $k$-algebra $A$ is both 
\textbf{left} and \textbf{right artinian}, that is, any descending chain of left and right ideals stops (the \textbf{DCC ascending chain condition}).

Thus, every Clifford algebra $\clpq{p}{q}$, as well as every group algebra $kG$, when $G$ is a finite group, which then makes $kG$ finite dimensional, have both chain conditions by a dimensionality argument. 
%In particular, this implies that every non-empty family of left ideals in these algebras contains a maximal element and a minimal element.

\begin{definition}
A left ideal $L$ in a ring $R$ is a \textbf{minimal left ideal} if $L\neq (0)$ and there is no left ideal $J$ with $(0) \subsetneq J \subsetneq L.$
\end{definition}

\noindent
One standard example of  minimal left ideals in matrix algebras $R=\Mat(n,k)$ are the subspaces 
$\mathrm{COL}(j)$, $1 \leq j \leq  n,$ of $\Mat(n,k)$ consisting of matrices $[a_{i,j}]$ such that $a_{i,k}=0$ when $k\neq j$ (cf. \cite[Example 7.9]{rotman}). 

The following proposition relates minimal left ideals in a ring $R$ to simple left $R$-modules. Recall that a left $R$-module $M$ is \textbf{simple} (or \textbf{irreducible}) if $M \neq \{0\}$ and $M$ has no proper nonzero submodules.

\begin{proposition}[Rotman \cite{rotman}]
\leavevmode
\begin{itemize}
\item[(i)] Every minimal left ideal $L$ in a ring $R$ is a simple left $R$-module.
\item[(ii)] If $R$ is left artinian, then every nonzero left ideal $I$ contains a minimal left ideal. 
\end{itemize}
\end{proposition}
Thus, the above proposition applies to Clifford algebras $\clpq{p}{q}$: every left spinor ideal $S$ in 
$\clpq{p}{q}$ is a simple left $\clpq{p}{q}$-module; and, every left ideal in $\clpq{p}{q}$ contains a spinor ideal.

Recall that if $D$ is a division ring, then a left (or right) $D$-module $V$ is called a 
\textbf{left} (or \textbf{right}) \textbf{vector space} over $D$. In particular, when the division ring is a field $k$, then we have a familiar concept of a $k$-vector space. Since the concept of linear independence of vectors generalizes from $k$-vector spaces to $D$-vector spaces, we have the following result.

\begin{proposition}[Rotman \cite{rotman}]
Let $V$ be a finitely generated\footnote{The term ``finitely generated" means that every vector in $V$ is a linear combination of a finite number of certain vectors $\{x_1,\ldots,x_n\}$ with coefficients from $R$. In particular, a $k$-vector space is finitely generated if and only if it is finite-dimensional \cite[Page 405]{rotman}.} left vector space over a division ring $D$. 
\begin{itemize}
\item[(i)] $V$ is a direct sum of copies of $D$; that is, every finitely generated left vector space over $D$ has a basis.
\item[(ii)] Any two bases of $V$ have the same number of elements.
\end{itemize}
\label{prop:propfg}
\end{proposition}
Since we know from the Structure Theorem, that every spinor left ideal~$S$ in simple Clifford algebras 
$\clpq{p}{q}$ ($p-q \neq 1 \bmod 4$) is a right $\BK$-module where $\BK$ is one of the division rings 
$\BR,\BC,$ or $\BH$, the above proposition simply tells us that every spinor left ideal $S$ is 
finite-dimensional over $\BK$ where $\BK$ is one of $\BR,\BC$ or $\BH.$ 

In semisimple Clifford algebras $\clpq{p}{q}$ ($p-q = 1 \bmod 4$), we have to be careful as the faithful double spinor representations are realized in the direct sum of two spinor ideals $S \oplus \hat{S}$ which are right $\BK \oplus \hat{\BK}$-modules, where $\BK =\BR$ or $\BH$.\footnote{Here, 
$\hat{S}= \{\hat{\psi} \mid \in S\}$, and similarly for $\hat{\BK}$, where $\hat{\phantom{S}}$ denotes the 
grade involution in $\clpq{p}{q}$.} Yet, it is easy to show that $\BK \oplus \hat{\BK}$ is not a division ring.

Thus, Proposition~\ref{prop:propfg} tells us that every finitely generated left (or right) vector space~$V$ over a division ring $D$ has a left (a right) dimension, which may be denoted $\dim V.$ In~\cite{jacobson} Jacobson gives an example of a division ring $D$ and an abelian group $V$, which is both a right and a left $D$-vector space, such that the left and the right dimensions are not equal. In our discussion, spinor minimal ideal $S$ will always be a left $\clpq{p}{q}$-module and a right $\BK$-module.

Since semisimple rings generalize the concept of a group algebra $\BC G$ for a finite group~$G$ (cf.~\cite{james,rotman}), we first discuss semisimple modules over a ring~$R$.   

\begin{definition}
A left $R$-module is \textbf{semisimple} if it is a direct sum of (possibly infinitely many) simple modules.
\end{definition}

The following result is an important characterization of semisimple modules.

\begin{proposition}[Rotman \cite{rotman}]
A left $R$-module $M$ over a ring $R$ is semisimple if and only if every submodule of $M$ is a direct summand.
\end{proposition}

Recall that if a ring $R$ is viewed as a left $R$-module, then its submodules are its left ideals, and, a left ideal is minimal if and only if it is a simple left $R$-module~\cite{rotman}.

\begin{definition}
A ring $R$ is \textbf{left semisimple}\footnote{One can define a \textbf{right semisimple} ring $R$ if it is a direct sum of minimal right ideals. However, it is known~\cite[Corollary 7.45]{rotman} that a ring is left semisimple if and only if it is right semisimple.} if it is a direct sum of minimal left ideals.
\end{definition}
\noindent

One of the important consequences of the above for the theory of Clifford algebras, is the following proposition.

\begin{proposition}[Rotman \cite{rotman}]
Let $R$ be a left semisimple ring.
\begin{itemize}
\item[(i)] $R$ is a direct sum of finitely many minimal left ideals.
\item[(ii)] $R$ has both chain conditions on left ideals.
\end{itemize}
\end{proposition}
From a proof of the above proposition one learns that,  while $R=\bigoplus_{i}L_i$, that is, $R$ is a direct sum of finitely-many left minimal ideals, the unity $1$ decomposes into a sum $1 = \sum_{i}f_i$ of mutually annihilating primitive idempotents $f_i$, that is, $(f_i)^2=f_i,$ and $f_if_j=f_jf_i=0, i\neq j.$ Furthermore, we find that $L_i=Rf_i$ for every $i.$

We can conclude from the following fundamental result~\cite{rotman,herstein} that every Clifford algebra $\clpq{p}{q}$ is a semisimple ring, because every Clifford algebra is a twisted group algebra $\BR^t[(\BZ_2)^n]$ for $n=p+q$ and a suitable twist~\cite{albuquerque,ablamowicz3,brasov}. 

\begin{theorem}[Maschke's Theorem]
If $G$ is a finite group and $k$ is a field whose characteristic does not divide $|G|$, the $kG$ is a left semisimple ring. 
\end{theorem}

For characterizations of left semisimple rings, we refer to~\cite[Section 7.3]{rotman}.

Before stating Wedderburn-Artin Theorem, which is all-important to the theory of Clifford algebras, 
we conclude this part with a definition and two propositions.

\begin{definition}
A ring $R$ is \textbf{simple} if it is not the zero ring and it has no proper nonzero two-sided ideals.
\end{definition}

\begin{proposition}[Rotman \cite{rotman}]
If $D$ is a division ring, then $R=\Mat(n,D)$ is a simple ring.
\end{proposition}

\begin{proposition}[Rotman \cite{rotman}]
If $R=\bigoplus_jL_j$ is a left semisimple ring, where the $L_j$ are minimal left ideals, then every simple $R$-module $S$ is isomorphic to $L_j$ for some $j.$
\label{prop:prop6}
\end{proposition}
\noindent
The main consequence of this last result is that every simple, hence irreducible, left $\clpq{p}{q}$-module, that is, every (left) spinor module of $\clpq{p}{q}$, is isomorphic to some minimal left ideal $L_j$ in the direct sum decomposition of $R=\clpq{p}{q}$.

Following Rotman, we divide the Wedderburn-Artin Theorem into the existence part and a uniqueness part. We also remark after Rotman that Wedderburn proved the existence theorem~\ref{thm:thm1} for semisimple $k$-algebras, where $k$ is a field, while E.~Artin generalized this result to what is now known as the Wedderburn-Artin Theorem.

\begin{theorem}[Wedderburn-Artin I]
A ring $R$ is left semisimple if and only if $R$ is isomorphic to a direct product of matrix rings over division rings $D_1,\ldots,D_m$, that is
\begin{gather}
R \cong \Mat(n_1,D_1) \times \cdots \times \Mat(n_m,D_m). \label{eq:eq1}
\end{gather}
\label{thm:thm1}
\end{theorem}
A proof of the above theorem yields that if $R=\bigoplus_jL_j$ as in Proposition~\ref{prop:prop6}, then each division ring $D_j=\End_R(L_j), j=1,\ldots,m,$ where $\End_R(L_j)$ denotes the ring of all $R$-endomorphisms  
of~$L_j$. Another consequence is the following corollary.

\begin{corollary}
A ring $R$ is left semisimple if and only if it is right semisimple.
\end{corollary}

Thus, we may refer to a ring as being \textbf{semisimple} without specifying from which side.\footnote{Not every simple ring is semisimple, cf.~\cite[Page 554]{rotman} and reference therein.} However, we have the following result which we know applies to Clifford algebras $\clpq{p}{q}.$ More importantly, its corollary explains part of the Structure Theorem which applies to simple Clifford algebras. Recall from the above that every Clifford algebra $\clpq{p}{q}$ is left artinian (because it is finite-dimensional).

\begin{proposition}[Rotman \cite{rotman}]
A simple left artinian ring $R$ is semisimple.
\end{proposition}

\begin{corollary}
If $A$ is a simple left artinian ring, then $A \cong \Mat(n,D)$ for some $n\geq 1$ and some division ring $D$.
\end{corollary}

Before we conclude this section with the second part of the Wedderburn-Artin Theorem, which gives certain uniqueness of the decomposition~(\ref{eq:eq1}), we state the following definition and a lemma.
\begin{definition}
Let $R$ be a left semisimple ring, and let 
\begin{gather}
R = L_1 \oplus \cdots \oplus L_n,
\end{gather}
where the $L_j$ are minimal left ideals. Let the ideals $L_1,\ldots,L_m$, possibly after re-indexing, be such that no two among them are isomorphic, and so that every $L_j$ in the given decomposition of $R$ is isomorphic to one and only one $L_i$ for $1 \leq i \leq m.$ The left ideals
\begin{gather}
B_i=\bigoplus_{L_j \cong L_i}L_j
\end{gather}
are called the \textbf{simple components} of $R$ relative to the decomposition $R=\bigoplus_jL_j.$
\end{definition}
\begin{lemma}[Rotman \cite{rotman}]
Let $R$ be a semisimple ring, and let 
\begin{gather}
R = L_1 \oplus \cdots \oplus L_n = B_1 \oplus \cdots \oplus B_m
\end{gather}
where the $L_j$ are minimal left ideals and the $B_i$ are the corresponding simple components of~$R$.
\begin{itemize}
\item[(i)] Each $B_i$ is a ring that is also a two-sided ideal in $R$, and $B_iB_j=(0)$ if $i \neq j.$
\item[(ii)] If $L$ is any minimal left ideal in $R$, not necessarily occurring in the given decomposition of $R$, then $L \cong L_i$ for some $i$ and $L\subseteq B_i.$
\item[(iii)] Every two-sided ideal in $R$ is a direct sum of simple components.
\item[(iv)] Each $B_i$ is a simple ring.  
\end{itemize}
\end{lemma}
\noindent
Thus, we will gather from the Structure Theorem, that for simple Clifford algebras $\clpq{p}{q}$ we have only one simple component, hence $m=1$, and thus all $2^k$ left minimal ideals generated by a complete set of $2^k$ primitive mutually annihilating idempotents which provide an orthogonal decomposition of the unity $1$ in $\clpq{p}{q}$ (see part (c) of the theorem and notation therein). Then, for semisimple Clifford algebras 
$\clpq{p}{q}$ we have obviously $m=2$.

Furthermore, we have the following corollary results.
\begin{corollary}[Rotman \cite{rotman}]
\leavevmode
\begin{itemize}
\item[(1)] The simple components $B_1,\ldots,B_m$ of a semisimple ring $R$ do not depend on a decomposition of $R$ as a direct sum of minimal left ideals; 
\item[(2)] Let $A$ be a simple artinian ring. Then, 
\begin{itemize}
\item[(i)] $A \cong \Mat(n,D)$ for some division ring $D$. If $L$ is a minimal left ideal in $A$, then every simple left $A$-module is isomorphic to $L$; moreover, $D^{\mathrm{op}} \cong \End_A(L).$\footnote{By 
$D^{\mathrm{op}}$ we mean the \textbf{opposite ring} of $D$: It is defined as 
$D^{\mathrm{op}}=\{a^{\mathrm{op}}\mid a \in D\}$ with multiplication defined as 
$a^{\mathrm{op}}\cdot b^{\mathrm{op}} = (ba)^{\mathrm{op}}.$}
\item[(ii)] Two finitely generated left $A$-modules $M$ and $N$ are isomorphic if and only if 
$\dim_D(M) =\dim_D(N).$ 
\end{itemize}
\end{itemize}
\end{corollary}
\noindent
As we can see, part (1) of this last corollary gives a certain invariance in the decomposition of a semisimple ring into a direct sum of simple components. Part (2i), for the left artinian Clifford algebras $\clpq{p}{q}$ implies that simple Clifford algebras ($p-q \neq 1 \bmod{4}$) are simple algebras isomorphic to a matrix algebra over a suitable division ring $D$. From the Structure Theorem we know that $D$ is one of 
$\BR,$ $ \BC$, or $\BH$, depending on the value of $p-q \bmod{8}.$ Part (2ii) tells us that any two spinor ideals $S$ and $S'$, which are simple right $\BK$-modules (due the right action of the division ring $\BK=f\clpq{p}{q}f$ on each of them) are isomorphic since their dimensions over $\BK$ are the same.

We conclude this introduction to the theory of semisimple rings with the following uniqueness theorem.
\begin{theorem}[Wedderburn-Artin II]
Every semisimple ring $R$ is a direct product,
\begin{gather}
R \cong \Mat(n_1,D_1) \times \cdots \times \Mat(n_m,D_m),
\label{thm:thm2}
\end{gather}
where $n_i \geq 1,$ and $D_i$ is a division ring, and the numbers $m$ and $n_i,$ as well as the division rings 
$D_i$, are uniquely determined by $R$.
\end{theorem}

Thus, the above results, and especially the Wedderburn-Artin Theorem (parts I and II), shed a new light on the main Structure Theorem given in the following section. In particular, we see it as a special case of the theory of semisimple rings, including the left artinian rings, applied to the finite dimensional Clifford algebras $\clpq{p}{q}$. 

We remark that the above theory applies to the group algebras $kG$ where $k$ is an algebraically closed field and $G$ is a finite group.

\section{The Main Structure Theorem on Real Clifford Algebras $\clpq{p}{q}$}
\label{sect:sect3}

We have the following main theorem that describes the structure of Clifford algebras $\clpq{p}{q}$ and their spinorial representations. In the following, we will analyze statements in that theorem. The same information is encoded in the well--known Table~1 in \cite[Page 217]{lounesto}. 

\begin{StructureTheorem}
Let $\cl_{p,q}$ be the universal real Clifford algebra over $(V,Q)$, $Q$ is non-degenerate of signature $(p,q)$.
\begin{itemize}
\item[(a)] When $p-q \neq 1 \bmod 4$ then $\cl_{p,q}$ is a simple algebra of dimension $2^{p+q}$ isomorphic with a full matrix algebra $\Mat(2^k, \BK)$ over a division ring $\BK$ where $k = q - r_{q-p}$ and $r_i$ is the Radon-Hurwitz number.\footnote{The Radon-Hurwitz number is defined by recursion as $r_{i+8}=r_i+4$ and these initial values: $r_0=0,$ $r_1=1,$ $r_2=r_3=2,r_4=r_5=r_6=r_7=3.$ \label{page:pageRH}} Here $\BK$ is one of $\BR, \BC$ or $\BH$ when $(p-q) \bmod 8$ is $0,2,$ or $3,7$, or $4,6$.
%%%%%%%%%%%%%%%%%%%%%%%%%%%%%%%%%%%%%%%%%%%%%%%%%%%%%%%
\item[(b)] When $p-q = 1 \bmod 4$ then $\cl_{p,q}$ is a semisimple algebra of dimension $2^{p+q}$ isomorphic to $\Mat(2^{k-1}, \BK) \oplus \Mat(2^{k-1}, \BK),$ $k = q - r_{q-p}$, and $\BK$ is isomorphic to $\BR$ or $\BH$ depending whether $(p-q) \bmod 8$ is $1$ or~$5$. Each of the two simple direct components of $\cl_{p,q}$ is projected out by one of the two central idempotents $\frac12(1\pm \be_{12 \ldots n}).$
%%%%%%%%%%%%%%%%%%%%%%%%%%%%%%%%%%%%%%%%%%%%%%%%%%%%%%%
\item[(c)] Any element $f$ in $\cl_{p,q}$ expressible as a product
\begin{equation}
f = \frac12(1\pm \be_{\iu_1})\frac12(1\pm \be_{\iu_2})\cdots\frac12(1\pm \be_{\iu_k})
\label{eq:f}
\end{equation}
where $\be_{\iu_j},$ $j=1,\ldots,k,$ are commuting basis monomials in $\cb{B}$ with square $1$ and $k = q - r_{q-p}$ generating a group of order $2^k$, is a primitive idempotent in $\cl_{p,q}.$ Furthermore, $\cl_{p,q}$ has a complete set of $2^k$ such primitive mutually annihilating idempotents which add up to the unity $1$ of $\cl_{p,q}$. 
%%%%%%%%%%%%%%%%%%%%%%%%%%%%%%%%%%%%%%%%%%%%%%%%%%%%%%%
\item[(d)] When $(p-q) \bmod 8$ is $0,1,2,$ or $3,7$, or $4,5,6$, then the division ring $\BK = f \cl_{p,q}f$ is isomorphic to $\BR$ or $\BC$ or $\BH$, and the map $S \times \BK \rightarrow S,$ $(\psi,\lambda) \mapsto \psi\lambda$ defines a right  $\BK$-module structure on the minimal left ideal $S=\cl_{p,q}f.$ 
%%%%%%%%%%%%%%%%%%%%%%%%%%%%%%%%%%%%%%%%%%%%%%%%%%%%%%%
\item[(e)] When $\cl_{p,q}$ is simple, then the map
\begin{equation}
\cl_{p,q} \stackrel{\gamma}{\longrightarrow} \End_\BK(S), \quad u \mapsto \gamma(u), \quad \gamma(u)\psi = u \psi
\label{eq:simple}
\end{equation}
gives an irreducible and faithful representation of $\cl_{p,q}$ in $S.$
%%%%%%%%%%%%%%%%%%%%%%%%%%%%%%%%%%%%%%%%%%%%%%%%%%%%%%%
\item[(f)] When $\cl_{p,q}$ is semisimple, then the map
\begin{equation}
\cl_{p,q} \stackrel{\gamma}{\longrightarrow} \End_{\BK \oplus \hat{\BK}}(S \oplus \hat{S}), \quad u \mapsto \gamma(u), \quad \gamma(u)\psi = u \psi
\label{eq:semisimple}
\end{equation}
gives a faithful but reducible representation of $\cl_{p,q}$ in the double spinor space $S \oplus \hat{S}$ where $S =\{ u f\, |\, u \in \cl_{p,q}\}$, $\hat{S} =\{ u \hat{f}\, |\, u \in \cl_{p,q}\}$ and 
$\hat{\phantom{m}}$ stands for the grade involution in $\cl_{p,q}.$ In this case, the ideal $S \oplus \hat{S}$ is a right $\BK \oplus \hat{\BK}$-module structure,  
$\hat{\BK} = \{\hat{\lambda} \,| \, \lambda \in \BK \}$, and $\BK \oplus \hat{\BK}$ is isomorphic to $\BR \oplus \BR$ when $p-q=1 \bmod 8$ or to $\BH \oplus \hat{\BH} $ when 
$p-q=5 \bmod 8.$
\end{itemize}
\label{th:structure}
\end{StructureTheorem}

Parts (a) and (b) address simple and semisimple Clifford algebras $\clpq{p}{q}$ which are distinguished by the value of $p-q \bmod{4}$ while the dimension of $\clpq{p}{q}$ is $2^{p+q}.$ 
For simple algebras, the Radon-Hurwitz number $r_i$ defined recursively as shown, determines the value of 
the exponent $k = q - r_{q-p}$ such that
\begin{gather}
\clpq{p}{q} \cong \Mat(2^k,\BK) \quad \text{when $p-q  \neq 1 \bmod{4}$}. 
\label{eq:simple}
\end{gather}
Then, the value of $p-q \bmod{8}$ (``Periodicity of Eight" cf.~\cite{lounesto,periodicity}) determines whether $\BK \cong \BR,\BC$ or~$\BH.$ Furthermore, this automatically tells us, based on the theory outlined above, that 
\begin{gather}
\clpq{p}{q} = L_1 \oplus \cdots \oplus L_N, \quad N=2^k,
\end{gather}
that is, that the Clifford algebras decomposes into a direct sum of $N=2^k$ minimal left ideals (simple left $\clpq{p}{q}$-modules) $L_i,$ each of which is generated by a primitive idempotent. How to find these primitive mutually annihilating idempotents, is determined in Part~(c).

In Part (b) we also learn that the Clifford algebra $\clpq{p}{q}$ is semisimple as it is the direct sum of two simple algebras:
\begin{gather}
\clpq{p}{q} \cong \Mat(2^{k-1},\BK) \oplus \Mat(2^{k-1},\BK) \quad \text{when $p-q  = 1 \bmod{4}$}. 
\label{eq:semisimple}
\end{gather}
Thus, we have two simple components in the algebra, each of which is a subalgebra. Notice that
the two algebra elements
\begin{gather}
c_1=\frac12(1 + \be_{12 \ldots n}) \quad \mbox{and} \quad c_2=\frac12(1 - \be_{12 \ldots n})
\label{eq:cidemps}
\end{gather}
are central, that is, each belongs to the center $Z(\clpq{p}{q})$ of the algebra.\footnote{The center 
$Z(A)$ of an $k$-algebra $A$ contains all elements in $A$ which commute with every element in~$A$. In particular, from the definition of the $k$-algebra, $\lambda 1 \in Z(\clpq{p}{q})$ for every 
$\lambda \in k$.} This requires that $n=p+q$ be odd, so that the unit pseudoscalar  $\be_{12 \ldots n}$ would commute with each generator $\be_i$, and that $(\be_{12 \ldots n})^2=1, $ so that 
expressions~(\ref{eq:cidemps}) would truly be idempotents. Notice, that the idempotents $c_1,c_2$ provide an 
\textbf{orthogonal decomposition} of the unity 1 since $c_1+c_2=1,$ and they are mutually annihilating since 
$c_1c_2=c_2c_1=0.$ Thus, 
\begin{gather}
\clpq{p}{q} = \clpq{p}{q}c_1 \oplus  \clpq{p}{q}c_2
\end{gather}
where each $\clpq{p}{q}c_i$ is a simple subalgebra of $\clpq{p}{q}.$ Hence, by Part (a), each subalgebra  is  
isomorphic to $\Mat(2^{k-1},\BK)$ where $\BK$ is either $\BR$ or $\BH$ depending on the value of $p-q \bmod{8},$ as indicated. 

Part (c) tells us how to find  a complete set of $2^k$ primitive mutually annihilating idempotents, obtained by independently varying signs $\pm$ in each factor in~(\ref{eq:f}), provide an orthogonal decomposition of the unity. The set of $k$ commuting basis monomials $\be_{\iu_1}, \ldots,\be_{\iu_k}$, which square to~$1,$ is not unique. Stabilizer groups of these $2^k$ primitive idempotents $f_1,\ldots,f_N$ $(N=2^k)$ under the conjugate action of Salingaros vee groups are discussed in~\cite{ablamowicz2,ablamowicz3}. It should be remarked, that each idempotent in~(\ref{eq:f}) must have exactly $k$ factors in order to be primitive.

Thus, we conclude from Part (c) that
\begin{gather}
\clpq{p}{q} = \clpq{p}{q}f_1 \oplus \cdots \oplus  \clpq{p}{q}f_N, \quad N=2^k,
\end{gather}
is a decomposition of the Clifford algebra $\clpq{p}{q}$ into a direct sum of minimal left ideals, or, simple left $\clpq{p}{q}$-modules.

Part (d) determines the unique division ring $\BK=f\clpq{p}{q}f$, where $f$ is any primitive idempotent, prescribed by the Wedderburn-Artin Theorem, such that the decomposition~(\ref{eq:simple}) 
or~(\ref{eq:semisimple}) is valid, depending whether the algebra is simple or not. This part also reminds us that the left spinor ideals, while remaining left $\clpq{p}{q}$ modules, are right $\BK$-modules. This is important when computing actual matrices in spinor representations (faithful and irreducible). Detailed computations of these representations in both simple and semisimple cases are shown in~\cite{spinorrepr}. Furthermore, package \texttt{CLIFFORD} has a built-in database which displays matrices representing generators of $\clpq{p}{q}$, namely $\be_1,\ldots,\be_n,$ $n=p+q$, for a certain choice of a primitive 
idempotent~$f$. Then, the matrix representing any element $u \in \clpq{p}{q}$ can the be found using the fact that the maps $\gamma$ shown on Parts (e) and (f), are algebra maps.

Finally, we should remark, that while for simple Clifford algebras the spinor minimal left ideal carries a \textbf{faithful} (and irreducible) representation, that is, $\ker \gamma = \{1\},$ in the case of semisimple algebras, each $\frac12$ spinor space $S$ and $\hat{S}$ carries an irreducible but not faithful representation. Only in the double spinor space $S \oplus \hat{S},$ one can realize the semisimple algebra faithfully. For all practical purposes, this means that each element $u$ in a semisimple algebra must be represented by a pair of matrices, according to the isomorphism~(\ref{eq:semisimple}). In practice, the two matrices can then be considered as a single matrix, but over $\BK\oplus\hat{\BK}$ which is isomorphic to $\BR \oplus \BR$ or $\BH \oplus \BH$, depending whether $p-q = 1 \bmod{8},$ or $p-q = 5 \bmod{8}.$ We have already remarked earlier that while $\BK$ is a division ring, $\BK\oplus\hat{\BK}$ is not.

\section{Conclusions} 
  
In this paper, the author has tried to show how the Structure Theorem on Clifford algebras $\clpq{p}{q}$ is related to the theory of semisimple rings, and, especially of left artinian rings. Detailed computations of spinor representations, which were distributed at the conference, came from~\cite{spinorrepr}.

\section{Acknowledgments}
Author of this paper is grateful to Dr. habil. Bertfried Fauser for his remarks and comments which have 
helped improve this paper.

\end{document}